\documentclass[twoside,12pt]{article}
\usepackage{a4, amsfonts, latexsym}
\begin{document}
\newtheorem{theorem}{Theorem}[section]
\newtheorem{prop}[theorem]{Proposition}
\newtheorem{lemma}[theorem]{Lemma}
\newtheorem{corollary}[theorem]{Corollary}
\newtheorem{conj}[theorem]{Conjecture}
\newcommand{\Q}{\mathbb Q}
\newcommand{\Z}{\mathbb Z}  
\newcommand{\R}{\mathbb R}
\newcommand{\F}{\mathbb F} 
\newcommand{\lb}{\left}
\newcommand{\rb}{\right}
\newcommand{\pf}{{\bf Proof: }}
\newcommand{\finpf}{{\hfill$\Box$}}
\newcommand{\mko}{\mathfrak{o}}
\newcommand{\mkp}{\mathfrak{p}}
\newcommand{\mkq}{\mathfrak{q}}
\newcommand{\mco}{\mathcal{O}}
\newcommand{\mcp}{\mathcal{P}}
\newcommand{\ul}{\underline}
\newcommand{\Gal}{\mathrm{Gal}} 
\newcommand{\tabeq}{\!\!\!\!&=&\!\!\!\!}
\newcommand{\tabneq}{\!\!\!\!&\neq&\!\!\!\!}
\newcommand{\tabplus}{\!\!\!\!&+&\!\!\!\!}
\newcommand{\tabminus}{\!\!\!\!&-&\!\!\!\!}
\newcommand{\tabspace}{\!\!\!\!& &\!\!\!\!}
\newcommand{\tabtimes}{\!\!\!\!&\times&\!\!\!\!}
\newcommand{\tabcolon}{\!\!\!\!&\colon&\!\!\!\!}
\newcommand{\tabdivide}{\!\!\!\!&\div&\!\!\!\!}
\newcommand{\tabimplies}{\!\!\!\!&\Rightarrow&\!\!\!\!}
\newcommand{\tabiff}{\!\!\!\!&\Leftrightarrow&\!\!\!\!}
\newcommand{\tabmapsto}{\!\!\!\!&\mapsto&\!\!\!\!}
\newcommand{\tabnotin}{\!\!\!\!&\notin&\!\!\!\!}
\begin{center}
\addtocounter{footnote}{1}
\footnotetext{1991 {\em Mathematics Subject Classification:} Primary 11G30; Secondary 11Y40.}
\addtocounter{footnote}{1}
\footnotetext{{\em Key words and phrases:} Curve of genus 2, Mordell-Weil rank, $\sqrt2$ multiplication, isogeny.} 
\large{\bf{THE MORDELL-WEIL RANK OF THE JACOBIAN OF A\\CURVE OF GENUS 2 WITH $\sqrt2$ MULTIPLICATION}}\\
\vspace{0.25in}
\normalsize{}
PETER R. BENDING\\
\vspace{0.25in}
\end{center}  
\begin{abstract}
We present a method for computing the Mordell-Weil rank of the jacobian of a curve of genus 2 with $\sqrt2$ multiplication, based on descent via isogenies of degree 2, and apply it to a family of curves.    
\end{abstract}             
\section{Introduction}
Let $A$ be an abelian variety defined over a perfect field $F$, and let $n$ be an integer greater than $1$; in this paper, $F$ will usually be a number field or a completion of a number field. In the case where $F$ is a number field $K$, the Mordell-Weil theorem states that $A(K)$ is finitely generated, the weak Mordell-Weil theorem states that $A(K)/nA(K)$ is finite, and it is easily shown that
\begin{equation}
\lb|\frac{A(K)}{nA(K)}\rb|=|A(K)[n]|\times n^r,\label{modn}
\end{equation}
where $r$ is the Mordell-Weil rank of $A$ over $K$. The following proposition, whose proof uses the theory of formal groups, is well-known:
\begin{prop}\label{formal}
Let $\mko$ be the ring of integers of $K$, and let $\mkp$ be a prime ideal of $\mko$. Suppose that $A$ has good reduction at $\mkp$, and that $n$ is a positive integer which is coprime to the characteristic of the residue field modulo $\mkp$. Then, the natural map from the group of points of order dividing $n$ on $A$ to the group of points on the reduction of $A$ modulo $\mkp$ is an embedding.   
\end{prop}
Due to this proposition, the order of $A(K)[n]$ can usually be easily computed, so the Mordell-Weil rank can usually be easily computed provided we know the order of $A(K)/nA(K)$. The Mordell-Weil rank, and so $A(K)/nA(K)$, is one of the main interests of computational algebraic number theory, motivating the study of the Selmer group, which is a group related to $A(K)/nA(K)$ and is effectively, but not necessarily practically, computable (involving local rather than global computations), and the Tate-Shafarevich group, which measures how close the Selmer group is to $A(K)/nA(K)$, and is much harder to compute. The latter is involved in the second Birch and Swinnerton-Dyer conjecture, which has been investigated widely for elliptic curves. 

We will now present an exact sequence which is of fundamental importance to descent methods for computing Mordell-Weil ranks. Let $\eta$ be an endomorphism of $A$ defined over $F$ whose kernel is finite (for example the multiplication by $n$ endomorphism, from now on denoted by $[n]$), $A'$ be another abelian variety defined over $F$ whose dimension is that of $A$, $\psi$ be an isogeny from $A$ to $A'$ defined over $F$ whose kernel is contained in that of $\eta$, and $\psi'$ be the isogeny from $A'$ to $A$ defined over $F$ such that $\psi'\circ\psi=\eta$. We will now explain how the computation of the order of $A(F)/\eta A(F)$ can be split into two computations, assuming that this group is finite. It is easily shown that the sequence of group homomorphisms
\begin{equation}
0\rightarrow\frac{A'(F)[\psi']}{\psi(A(F)[\eta])}\rightarrow\frac{A'(F)}{\psi A(F)}\stackrel{\psi'}{\rightarrow}\frac{A(F)}{\eta A(F)}\rightarrow\frac{A(F)}{\psi'A'(F)}\rightarrow0\label{seq}
\end{equation}
is exact. The order of the second group can usually be easily computed, so the order of $A(F)/\eta A(F)$ can usually be easily computed provided we know the orders of $A'(F)/\psi A(F)$ and $A(F)/\psi'A'(F)$ (the latter two groups are finite since the former group is).

The idea described above of splitting one computation into two has been applied in the case where $F$ is a number field $K$, to compute Mordell-Weil ranks over $K$, as follows. For elliptic curves, it has been applied with $\eta$ being $[2]$ and $\psi,\psi'$ being isogenies of degree $2$ (see~\cite{cremona}, pp. 63-68). For jacobians of curves of genus $2$, it has been applied with $\eta$ being $[2]$ and $\psi,\psi'$ being isogenies of degree $4$, called Richelot isogenies (see~\cite{CF}, Chapters 10 and 11). In this paper, for jacobians of curves of genus $2$ with $\sqrt2$ multiplications, it will be applied twice: once with $\eta$ being $[2]$ and $\psi,\psi'$ being a $\sqrt2$ multiplication $\varepsilon$ (with $A'=A$), and once with $\eta$ being $\varepsilon$ and $\psi,\psi'$ being isogenies of degree $2$.

The method presented in this paper is relevant to the well-known conjecture that
 the jacobian of a curve of genus $2$ defined over $\Q$ which is simple over $\Q
$ and of $GL_2$-type is isogenous over $\Q$ to a factor of the jacobian of a modular curve. 
\bigskip\\ 
I would like to thank the EPSRC for supporting this research, the University of Kent at Canterbury for its hospitality, and Dr. J. R. Merriman for his interest and encouragement. 

\section{The equation giving the Mordell-Weil rank}\label{secprop}
We will now establish notation which will be used throughout this paper. Let $C,C'$ be curves defined over a perfect field $F$ whose jacobians $J,J'$ have $\sqrt2$ multiplications $\varepsilon,\varepsilon'$ defined over $F$ killing points $P_0,P_0'$ respectively of order $2$ defined over $F$, and suppose that there are isogenies $\phi,\phi'$ defined over $F$ from $J$ to $J'$, $J'$ to $J$ whose kernels are generated by $P_0,P_0'$ respectively, with the property that $\phi'\circ\phi=\varepsilon$.

There are group embeddings
\begin{equation}
\lambda\colon\frac{J'(F)}{\phi J(F)}\hookrightarrow\frac{F^*}{{F^*}^2},\lambda'\colon\frac{J(F)}{\phi'J'(F)}\hookrightarrow\frac{F^*}{{F^*}^2},\label{embed}
\end{equation}
defined as follows. Let $\overline{F}$ be an algebraic closure of $F$. We define $\lambda$ to be the embedding which sends the point $P'$ in $J'(F)$ to the element $[D]$ with the property that $F(\sqrt{D})$ is the field of definition over $F$ of the two points $P_1,P_2$ in $J(\overline{F})$ such that $\phi(P_1)=\phi(P_2)=P'$. Similarly, we define $\lambda'$ to be the embedding which sends the point $P$ in $J(F)$ to the element $[D]$ with the property that $F(\sqrt{D})$ is the field of definition over $F$ of the two points $P_1',P_2'$ in $J'(\overline{F})$ such that $\phi'(P_1')=\phi'(P_2')=P$. Let $I,I'$ be the images of $\lambda,\lambda'$ respectively, isomorphic as groups to $J'(F)/\phi J(F),J(F)/\phi'J'(F)$ respectively.

For the case where $F$ is a number field $K$, there is another way of describing the embeddings above. Let $\overline{K}$ be an algebraic closure of $K$. The coboundary map from $J'(K)$ to $H^1(\Gal(\overline{K}/K),J[\phi])$ induces an embedding from $J'(K)/\phi J(K)$ to $H^1(\Gal(\overline{K}/K),J[\phi])$, and the latter group is isomorphic to $K^*/{K^*}^2$ by Hilbert 90 (since $J[\phi]$ is isomorphic as $\Gal(\overline{K}/K)$-modules to $\mu_2$), so we get an embedding from $J'(K)/\phi J(K)$ to $K^*/{K^*}^2$. Similarly, we get an embedding from $J(K)/\phi'J'(K)$ to $K^*/{K^*}^2$. These embeddings can be generalised for isogenies of prime power degree, as discussed in~\cite{schaefer}. The fact that the degrees of $\phi$ and $\phi'$ are $2$ means that deriving explicit conditions for the existence of a particular element in $I$ or in $I'$ (as done in the next section) is relatively easy.   
         
For the case where $F$ is a number field $K$, we will state and prove a theorem which relates the Mordell-Weil rank of $J$ and $J'$ over $K$ to the product of the orders of $I$ and $I'$ (since $J$ and $J'$ are isogenous over $K$, their ranks are the same). The following lemma, for perfect $F$ in general, allows us to do this. Because the RHS is always $2^4$, we do not have to compute the terms in the LHS when working out Mordell-Weil ranks.   
\begin{lemma}\label{initlem}
\[|J(F)[2]|\lb|\frac{J'(F)[\phi']}{\phi(J(F)[\varepsilon])}\rb|^2\lb|\frac{J(F)[\varepsilon]}{\varepsilon(J(F)[2])}\rb|=2^4.\]
\end{lemma}
\pf The orders of the groups $J(F)[2]$ and $J(F)[\varepsilon]$ determine the orders of the groups $J'(F)[\phi']/\phi(J(F)[\varepsilon])$ and $J(F)[\varepsilon]/\varepsilon(J(F)[2])$, as shown in the following table, which lists all possible orders of the former two groups. In all the rows, the product of the first number, the square of the third number and the fourth number is $2^4$, establishing the lemma.

Some points can be made about the possible orders of the former two groups: $2$ divides the order of both of them since both contain the point $P_0$ of order $2$, the second order divides the first order since the second group is contained in the first group, and it is impossible for the first order to be $2^1$ and the second order to be divisible by $2^3$ simultaneously (if this was true, then $P_0$ would be the only non-trivial point in $J(F)$ killed by $\varepsilon$, and all of the points of order $2$ in $J(F)$ not killed by $\varepsilon$, six or fourteen in total, would be sent by $\varepsilon$ to $P_0$, a contradiction since the degree of $\varepsilon$ is $4$).
\begin{center}
\begin{tabular}{ccccc}\\
$|J(F)[2]|$&$|J(F)[\varepsilon]|$&&$\lb|\frac{J'(F)[\phi']}{\phi(J(F)[\varepsilon])}\rb|$&$\lb|\frac{J(F)[\varepsilon]}{\varepsilon(J(F)[2])}\rb|$\\\hline
$2^1$&$2^1$&&$2^1$&$2^1$\\
$2^2$&$2^1$&&$2^1$&$2^0$\\
$2^2$&$2^2$&&$2^0$&$2^2$\\
$2^3$&$2^2$&&$2^0$&$2^1$\\
$2^4$&$2^2$&&$2^0$&$2^0$
\end{tabular}
\end{center}
\finpf
\bigskip\\
For the case where $F$ is a number field $K$, the following theorem relates the Mordell-Weil rank of $J$ and $J'$ over $K$ to the product of the orders of $I$ and $I'$:   
\begin{theorem}\label{prop}
$I$ and $I'$ are finite, and 
\begin{equation} 
2^{4+r}=|I|^2|I'|^2,\label{result}
\end{equation}  
where $r$ is the Mordell-Weil rank of $J$ and $J'$.  

Let $S,S'$ be subsets of $I,I'$ whose elements generate $I,I'$ respectively. For each element $[D]$ of $S$, let $P'_{[D]}$ be a point on $J'$ defined over $K$ such that $\lambda(P'_{[D]})=[D]$. Similarly, for each element $[D]$ of $S'$, let $P_{[D]}$ be a point on $J$ defined over $K$ such that $\lambda'(P_{[D]})=[D]$. Then, the points $\phi'(P'_{[D]}),\varepsilon\circ\phi'(P'_{[D]})$ ($[D]$ ranging across the elements of $S$) and the points $P_{[D]},\varepsilon(P_{[D]})$ ($[D]$ ranging across the elements of $S'$) generate a subgroup of odd index of the Mordell-Weil group of $J$.  
\end{theorem}
\pf From (\ref{modn}) with $n=2$ and by applying the exact sequence (\ref{seq}) twice, firstly with $A,A',\psi,\psi',\eta$ being $J,J,\varepsilon,\varepsilon,[2]$ respectively, and secondly with $A,A',\psi,\psi',\eta$ being $J,J',\phi,\phi',\varepsilon$ respectively, we see that $J'(K)/\phi J(K),J(K)/\phi'J'(K)$ are finite (since $A(K)/2A(K)$ is), and that 
\[|J(K)[2]|\times2^r=\frac{\lb|\frac{J'(K)}{\phi J(K)}\rb|^2\lb|\frac{J(K)}{\phi'J'(K)}\rb|^2}{\lb|\frac{J'(K)[\phi']}{\phi(J(K)[\varepsilon])}\rb|^2\lb|\frac{J(K)[\varepsilon]}{\varepsilon(J(K)[2])}\rb|}.\]
$I,I'$ are isomorphic as groups to $J'(K)/\phi J(K),J(K)/\phi'J'(K)$ respectively, so (\ref{result}) now follows from Lemma~\ref{initlem}. By applying the exact sequence (\ref{seq}) twice in the same way, we deduce the rest of the theorem. 
\finpf
\bigskip\\
To apply Theorem~\ref{prop}, we will need methods of working out $I$ and $I'$, and explicit examples of curves $C$ and $C'$. These topics will be discussed in the next two sections. For application in the next two sections, we note that $C,C'$ have Weierstrass models
\[C\colon Y^2=F_2(X)F_4(X),C'\colon Y^2=F_2'(X)F_4'(X),\]
where $F_2,F_2',F_4,F_4'$ are polynomials over $\mko$, the ring of integers of $K$, of degrees less than or equal to $2,2,4,4$ respectively, with the points $P_0,P_0'$ on $J,J'$ corresponding to the zeros of $F_2,F_2'$ respectively.

\section{The groups $I$ and $I'$}\label{groups} 
In this section, we discuss methods of working out $I$ and $I'$. This task may not be easy in general; in this paper, we will concentrate on curves $C$ and $C'$ such that $I$ and  $I'$ are generated mostly by elements $[D]$ ($D$ an element of $K^*$) with at least one of the following properties:
\begin{itemize}
\item $C(K)$ (resp. $C'(K)$) is non-empty, and the zeros of $F_2$ (resp. $F_2'$) are defined over $K(\sqrt{D})$ individually.
\item The zeros of $F_4$ (resp. $F_4'$) can be partitioned into two pairs of two points, the pairs being defined over $K(\sqrt{D})$ but not over $K$, and conjugate over $K$. 
\end{itemize}
(Note that, if an element $[D]$ has at least one of these properties, it is clearly in $I$ (resp. $I'$)).
\bigskip\\
We will now state and prove five lemmas which are useful for narrowing down $I$ and $I'$. The first two lemmas give explicit conditions on $F_2,F_4$ (resp. $F_2',F_4'$) for the existence of an element $[D]$ in $I$ (resp. $I'$), the third lemma discusses the reductions of $J$ and $J'$ at the prime ideals of $K$ given the existence of an element $[D]$ in $I$ or $I'$, the fourth lemma compares the embeddings (\ref{embed}) for $K$ with those for a completion of $K$ w.r.t. a valuation associated to a prime ideal of $K$, and the fifth lemma compares the embeddings (\ref{embed}) for $K$ with those for $\R$.

First, we will establish standard notation from algebraic number theory. Letting $\mkp$ be a prime ideal of $\mko$, we define:
\begin{itemize}
\item $v_\mkp$ to be the valuation of $K$ associated to $\mkp$;
\item $\mco_\mkp$ to be the valuation ring of $v_\mkp$;
\item $\mcp_\mkp$ to be the valuation ideal of $v_\mkp$;
\item $K^c_\mkp$ to be a completion of $K$ w.r.t. $v_\mkp$;
\item $v^c_\mkp$ to be the valuation of $K^c_\mkp$ extending $v_\mkp$;
\item $\mco^c_\mkp$ to be the valuation ring of $v^c_\mkp$;
\item $\mcp^c_\mkp$ to be the valuation ideal of $v^c_\mkp$.
\end{itemize}   
\begin{lemma}\label{lem2}
Suppose that $[D]$ is an element of $I$ (resp. $I'$). Then, for all prime ideals $\mkp$ of $\mko$ such that $v_\mkp(D)$ is odd (note that this property is independent of the choice of $D$), $F_2$ or $F_4$ (resp. $F_2'$ or $F_4'$) is a constant times a square modulo $\mkp$.
\end{lemma}
\pf By multiplying $D$ by an appropriate square in $K$, we may assume that $v_\mkp(D)=1$, i.e. that $D$ is a uniformising element of $\mcp_\mkp$. The proofs of this lemma for $I,F_2,F_4$ and for $I',F_2',F_4'$ are analogous, so we will concentrate on the former proof.

There are two points $P_1,P_2$ in $J(K(\sqrt{D}))$, not defined but conjugate over $K$, such that $\phi(P_1)=\phi(P_2)$, which implies that
\begin{equation}
P_2=P_0+P_1\label{eqn}  
\end{equation}
(recall that the kernel of $\phi$ is generated by $P_0$). 

Write $P_i=\underline{x_{i1}}-\underline{x_{i2}}$, where $\underline{x_{ij}}$ is a point in $C(\overline{K})$ ($\overline{K}$ is an algebraic closure of $K$). There are two cases to consider:
\begin{enumerate}
\item Let $\iota$ be the hyperelliptic involution on $C$. Two of the points
\begin{equation}\ul{x_{01}},\iota(\ul{x_{02}}),\ul{x_{11}},\iota(\ul{x_{12}}),\iota(\ul{x_{21}}),\ul{x_{22}}\label{list}
\end{equation}
are sent to each other by $\iota$.
\item The above doesn't happen.
\end{enumerate}
The second case is the more usual case, and the more complicated case as it involves a function on $C$.
\bigskip\\
Consider the first case. Now none of the $P_i$'s are zero, since the kernel of $\phi$ is generated by $P_0$ and $P_1,P_2$ are not defined over $K$, so the first two points, the middle two points, and the last two points in (\ref{list}) are
not sent to each other. If any other two of the points in (\ref{list}) are sent
to each other, then, as we will explain below, the ratio of $D$ and the discriminant $\Delta$ of $F_2$ lies in ${K^*}^2$, which implies that $\Delta$ is an element of $\mkp$, since it is an element of $\mko$ and $D$ is a uniformising element of $\mcp_\mkp$, which implies that $F_2$ is a constant times a square modulo $\mkp$.

Suppose that any other two of the points in (\ref{list}) are sent to each other; we will concentrate on the points $\ul{x_{01}},\ul{x_{11}}$, the other pairs being treated similarly. For this pair, $P_2=[\iota(\ul{x_{02}})-\ul{x_{12}}]$ from (\ref{eqn}), so, if an automorphism $\sigma$ of $\overline{K}$ over $K$ varies $\sqrt{D}$, then, since $P_1,P_2$ are defined over $K(\sqrt{D})$ but not over $K$, and conjugate
over $K$, we see that $\sigma(\iota(\ul{x_{02}}))=\ul{x_{11}},\sigma(\ul{x_{12}})=\ul{x_{12}}$ or $\sigma(\iota(\ul{x_{02}}))=\iota(\ul{x_{12}}),\sigma(\ul{x_{12}})=\iota(\ul{x_{11}})$. If the second possibility happens, then $\sigma^2(\ul{x_{02}})=\ul{x_{01}}$ (dropping the $\iota$'s since $\ul{x_{01}},\ul{x_{02}}$ are Weierstrass points on $C$), which implies that $\ul{x_{02}}=\ul{x_{01}}$ since $\ul{x_{02}}$ is defined over at worst a quadratic extension of $K$, which is absurd. If the first possibility happens, then $\sigma(\ul{x_{02}})=\ul{x_{01}}$
(again dropping the $\iota$'s). This holds for all automorphisms of $\overline{K}$ over $K$ varying $\sqrt{D}$, so the ratio of $D$ and the discriminant of $F_2$ lies in ${K^*}^2$, as required.
\bigskip\\
Consider the second case. By (\ref{eqn}), the hypotheses of this case, and the fact that $P_1,P_2$ are defined over $K(\sqrt{D})$ but not over $K$, and conjugate over $K$, the divisor $\ul{x_{01}}-\ul{x_{02}}+\ul{x_{11}}-\ul{x_{12}}-\ul{x_{21}}+\ul{x_{22}}$ is the divisor of the function
\[(Y-\sqrt{D}(P+QX+RX^2+SX^3))/((X-x_{02})(X-x_{12})(X-x_{21}))\]
for some elements $P,Q,R,S$ of $K$, defining the element $x_{ij}$ of $\overline{K}\cup\infty$ to be the $X$-coordinate of $\ul{x_{ij}}$. So we have
\[F_2(X)F_4(X)-D(P+QX+RX^2+SX^3)^2=N\prod_{i=0}^2\prod_{j=0}^1(X-x_{ij}),\]
for some element $N$ of $K^*$, which implies that $P+QX+RX^2+SX^3=F_2(X)(UX+V)$
for some elements $U,V$ of $K$, and that
\begin{equation}
F_4(X)-DF_2(X)(UX+V)^2=N\prod_{i=1}^2\prod_{j=0}^1(X-x_{ij}).\label{ident}
\end{equation}
Write $U=\frac{U'}{T},V=\frac{V'}{T}$, where $U',V',T$ are elements of $\mco_\mkp$ not all in $\mcp_\mkp$. By (\ref{ident}), we have
\begin{equation}
T^2F_4(X)-DF_2(X)(U'X+V')^2=T^2N\prod_{i=1}^2\prod_{j=0}^1(X-x_{ij}).\label{ident2}
\end{equation}
This equation is the key to the proof for the second case.

Suppose that $T$ is not an element of $\mcp_\mkp$. We will show that $F_4$ is a
constant times a square modulo $\mcp_\mkp$; then, since $F_4$ is a polynomial over $\mko$, it will be a constant times a square modulo $\mkp$, as required. Multiplying both sides of (\ref{ident2}) by $T^{-2}$, an element of $\mco_\mkp$ since $T$ is not an element of $\mcp_\mkp$, we see that it suffices to show that
\begin{equation}
N\prod_{i=1}^2\prod_{j=0}^1(X-x_{ij}),\label{poly}
\end{equation}
a polynomial over $\mco_\mkp$, congruent modulo $\mcp_\mkp$ to $F_4$ since $D$ is an element of $\mcp_\mkp$, is a constant times a square modulo $\mcp_\mkp$. Since $P_1,P_2$ are defined over $K(\sqrt{D})$ but not over $K$, and conjugate over $K$, (\ref{poly}) multiplied by some unit of $\mco_\mkp$ is
\begin{equation}
\prod_{s\in\{\pm1\}}(\alpha_1+s\alpha_2\sqrt{D}+X(\beta_1+s\beta_2\sqrt{D})+X^2(\gamma_1+s\gamma_2\sqrt{D}))\label{poly2}
\end{equation}
for some elements $\alpha_1,\alpha_2,\beta_1,\beta_2,\gamma_1,\gamma_2$ of $K$,
a polynomial over $\mco_\mkp$; it suffices to show that (\ref{poly2}) is a constant times a square modulo $\mcp_\mkp$. Taking into account the fact that $D$ is
a uniformising element of $\mcp_\mkp$, we see that $\alpha_1,\alpha_2,\beta_1,\beta_2,\gamma_1,\gamma_2$ are all elements of $\mco_\mkp$; then, (\ref{poly2}) is congruent to the square of $(\alpha_1+\beta_1X+\gamma_1X^2)$ modulo $\mcp_\mkp$.

Suppose that $T$ is an element of $\mcp_\mkp$. We will show that $F_2$ is a constant times a square modulo $\mcp_\mkp$; then, since $F_2$ is a polynomial over $\mko$, it will be a constant times a square modulo $\mkp$, as required. It suffices to show that $(U'X+V')^2F_2(X)$ is a constant times a square modulo $\mcp_\mkp$; indeed, $U',V'$ are not both elements of $\mcp_\mkp$, since $T$ is an element and $U',V',T$ are not all elements. Dividing both sides of (\ref{ident2}) by
$D$, we see that it suffices to show that
\begin{equation}
T^2ND^{-1}\prod_{i=1}^2\prod_{j=0}^1(X-x_{ij}),\label{poly3}
\end{equation}
a polynomial over $\mco_\mkp$ (note that $T^2D^{-1}$ is an element of $\mcp_\mkp$, and so of $\mco_\mkp$, since $T$ is an element and $D$ is a uniformising element), congruent modulo $\mcp_\mkp$ to $(U'X+V')^2F_2(X)$ since $T^2D^{-1}$ is an element of $\mcp_\mkp$, is a constant times a square modulo $\mcp_\mkp$. The proof of this for (\ref{poly3}) is similar to the proof of this for (\ref{poly2}), which was given above.
\finpf
\begin{lemma}\label{lem3}
Suppose that there is an embedding $\rho\colon K\hookrightarrow\R$, and that $[D]$ is an element of $I$ (resp. $I'$) with the property that $\rho(D)$ is negative (note that this property is independent of the choice of $D$). Then, it is impossible for the leading, constant coefficients of $\rho(F_2)$ (resp. $\rho(F_2')$) to be non-zero and have different signs, the leading, constant coefficients of $\rho(F_4)$ (resp. $\rho(F_4')$) to be non-zero and have different signs, and the
leading coefficients of $\rho(F_2),\rho(F_4)$ (resp. $\rho(F_2'),\rho(F_4')$) to have the same sign.
\end{lemma}
\pf The proof of this lemma will be based on the proof of Lemma~\ref{lem2}. As in Lemma~\ref{lem2}, the proofs of this lemma for $I,F_2,F_4$ and for $I',F_2',F_4'$ are analogous, so we will concentrate on the former proof.

Suppose, for a contradiction, that the alleged impossibility is in fact possible. As in Lemma~\ref{lem2}, there are two cases to consider: the case where two of the points in (\ref{list}) are sent to each
other by $\iota$, and the case where this doesn't happen.
\bigskip\\
Consider the first case. As in Lemma~\ref{lem2}, the ratio of $D$ and the discriminant $\Delta$ of $F_2$ lies in ${K^*}^2$. So $\rho(\Delta)$ is negative, since $\rho(D)$ is. But, by our supposition, the leading, constant coefficients of $\rho(F_2)$ are non-zero and have different signs, and so the discriminant of $\rho(F_2)$, which is $\rho(\Delta)$, is positive. So the first case cannot occur.
\bigskip\\
Consider the second case. As in Lemma~\ref{lem2}, (\ref{ident}) holds. Attacking the coefficients of both sides of (\ref{ident}) with $\rho$, we see that our supposition is impossible, taking into account the fact that $\rho(D)$ is negative.
\finpf
\begin{lemma}\label{lem1}
Suppose that $[D]$ is an element of $I\cup I'$. Then, for all prime ideals $\mkp$ of $\mko$ not containing $2$ such that $K(\sqrt{D})/K$ is ramified at $\mkp$, $J$ and $J'$ have bad reduction at $\mkp$.
\end{lemma}
\pf The case where $[D]$ is an element of $I'$ can be proved in an analogous way to
the case where $[D]$ is an element of $I$, so we will concentrate on the former. Since $J$ and $J'$ are isogenous over $K$, their primes of bad reduction are the same, so it suffices to show that $J$ has bad reduction at $\mkp$.  

Since $K(\sqrt{D})/K$ is ramified at $\mkp$, $K^c_\mkp(\sqrt{D})/K^c_\mkp$ is ramified at $\mcp^c_\mkp$. Let $w^c_\mkp$ be the valuation of $K^c_\mkp(\sqrt{D})$ extending $v^c_\mkp$. Then, by Theorems 24 and 23 in~\cite{FT}, there is an element $\pi$ of $K^c_\mkp(\sqrt{D})$ such that the valuation ring, ideal of $w^c_\mkp$ is $\mco^c_\mkp[\pi],(\mcp^c_\mkp,\pi)$ respectively.   
%

Suppose, for a contradiction, that $J$ has good reduction at $\mkp$. Then, there is a way of embedding $J$ over $K^c_\mkp$ into projective $N$-space over $K^c_\mkp$ for some positive integer $N$, so that the reduction of the image is an abelian variety defined over $\mco^c_\mkp/\mcp^c_\mkp$.

There are two points $P_1,P_2$ in $J(K^c_\mkp(\sqrt{D}))$, not defined but conjugate over $K^c_\mkp$, such that $\phi(P_1)=\phi(P_2)$. Since the valuation ring of $w^c_\mkp$ is $\mco^c_\mkp[\pi]$, we may write
\[P_1=(a_0+b_0\pi:\cdots:a_N+b_N\pi),P_2=(a_0+b_0\sigma(\pi):\cdots:a_N+b_N\sigma(\pi)),\] 
where the $a_i$'s and the $b_i$'s are in $\mco^c_\mkp$, and $\sigma$ is the non-trivial automorphism of $K^c_\mkp(\sqrt{D})$ over $K^c_\mkp$. Since the valuation ideal of $w^c_\mkp$ is $(\mcp^c_\mkp,\pi)$, $P_1,P_2$ are the same modulo $(\mcp^c_\mkp,\pi)$, which contradicts the fact that $J$ has good reduction at $\mkp$, by (\ref{eqn}) and Proposition~\ref{formal}; indeed, since $(\mcp^c_\mkp,\pi)$ does not contain $2$, $P_0$ (being a torsion point of order $2$) and the identity point are not the same modulo $(\mcp^c_\mkp,\pi)$.
\finpf
\begin{lemma}\label{lem4}
Let $\mkp$ be a prime ideal of $\mko$. Then, the product of the orders of the images of $I$ and $I'$ under the natural map $K^*/{K^*}^2\rightarrow{K^c_\mkp}^*/{{K^c_\mkp}^*}^2$ is less than or equal to $2^2$ if $\mkp$ does not contain $2$, and less than or equal to $2^{2+[K^c_\mkp\colon\Q_2]}$ otherwise.
\end{lemma}
\pf By Proposition 2.4 in~\cite{schaefer}, the order of $J(K^c_\mkp)/2J(K^c_\mkp)$ is $\Phi$ times the order of $J(K^c_\mkp)[2]$, where $\Phi=1$ if $\mkp$ does not contain $2$, and $\Phi=2^{2[K^c_\mkp\colon\Q_2]}$ otherwise. 

By applying the exact sequence (\ref{seq}) twice, exactly as in Theorem~\ref{prop}, we see that $J'(K^c_\mkp)/\phi J(K^c_\mkp),J(K^c_\mkp)/\phi'J'(K^c_\mkp)$ are finite (since $J(K^c_\mkp)/2J(K^c_\mkp)$ is), and that 
\[\Phi|J(K^c_\mkp)[2]|=\frac{\lb|\frac{J'(K^c_\mkp)}{\phi J(K^c_\mkp)}\rb|^2\lb|\frac{J(K^c_\mkp)}{\phi'J'(K^c_\mkp)}\rb|^2}{\lb|\frac{J'(K^c_\mkp)[\phi']}{\phi(J(K^c_\mkp)[\varepsilon])}\rb|^2\lb|\frac{J(K^c_\mkp)[\varepsilon]}{\varepsilon(J(K^c_\mkp)[2])}\rb|}.\]
Since the images of $I,I'$ under the natural map $K^*/{K^*}^2\rightarrow{K^c_\mkp}^*/{{K^c_\mkp}^*}^2$ embed into $J'(K^c_\mkp)/\phi J(K^c_\mkp),J(K^c_\mkp)/\phi'J'(K^c_\mkp)$ respectively, the lemma follows from Lemma~\ref{initlem}.
\finpf
\begin{lemma}\label{lem5}
Suppose that there is an embedding $\rho\colon K\hookrightarrow\R$. Then, the product of the orders of the images of $I$ and  $I'$ under the map $K^*/{K^*}^2\rightarrow\R^*/{\R^*}^2$ induced by $\rho$ is less than or equal to $2$.
\end{lemma}
\pf The proof of this lemma will be based on the proof of Lemma~\ref{lem4}. Throughout the proof, we regard $J$ and $J'$ as being defined over $\R$ via $\rho$. By Proposition 2.5 in~\cite{schaefer}, the order of $J(\R)/2J(\R)$ is the order of $J(\R)[2]$ divided by $4$. As in Lemma~\ref{lem4}, the lemma now follows from Lemma~\ref{initlem}.
\finpf

\section{Explicit examples}
We have been able to prove the proposition and the lemmas in the previous two sections without explicit examples of curves $C,C'$, but to apply these results we will need explicit examples. The following theorem, proved and discussed in Chapter 4 of~\cite{bending}, provides such examples:
\begin{prop}\label{isogthm}
Let $U,V,W,\Delta\in K$, and let $F_6,F_6'$ be the polynomials defined by
\begin{eqnarray*}
F_6(X)\colon\tabeq\Delta UV(X^2+UX+V)\\
\tabtimes\lb[UVX^4+V\lb(\frac{1}{4}W(U-V)(U+V)+V^2+4\rb)X^3\rb.\\
\tabplus U\lb(\frac{1}{4}W(U-V)(U+V)+V^2+VW-4\rb)X^2\\
\tabplus(W(U^2+V^2)-4(V^2+4))X+UV(W-4)],\\
F_6'(X)\colon\tabeq\Delta UV(V-2)((U-2)X^2+2(V-U+2)X+2(U-V))\\
\tabtimes[(W(V-1)(U-V)^2-4(V(U-V)^2-(V-2)^2))X^4\\
\tabplus2(W(U-V)(2(1-V)U+V(3V-2))\\
\tabplus4(2UV(U-2(V+1))+(3V-2)(V^2+4)))X^3\\
\tabplus(W(U((V^2+4(V-1))U-6V(3V-2))-V^2(V^2-2(7V-2)))\\
\tabplus4(4UV(3(V+2)-2U)+(V^2+4)(V^2-2(7V-2))))X^2\\
\tabminus2V(W(U(UV+4(1-2V))+V^2(6-V))\\
\tabminus4(4U(U-V-4)-(V-6)(V^2+4)))X\\
\tabplus WV^2(U-V)(U+V-4)-4V(4U(U-4)-(V-4)(V^2+4))].
\end{eqnarray*}
Suppose that $F_6,F_6'$ are non-trivial and have no multiple zeros, and let $C,C'$ be the curves of genus $2$ defined by
\[C\colon Y^2=F_6(X),C'\colon Y^2=F_6'(X).\]
Then, the jacobians $J,J'$ of $C,C'$ have $\sqrt2$ multiplications $\varepsilon,\varepsilon'$ defined over $K$ killing the points $P_0,P_0'$ of order $2$ defined over $K$ corresponding to the zeros of
\[X^2+UX+V,(U-2)X^2+2(V-U+2)X+2(U-V)\]
respectively, and there are isogenies $\phi,\phi'$ defined over $K$ from $J$ to $J'$, $J'$ to $J$ whose kernels are generated by $P_0,P_0'$ respectively, with the property that $\phi'\circ\phi=\varepsilon$ (so the curves $C,C'$ satisfy the conditions described at the beginning of section~\ref{secprop}, with $F=K$).
\end{prop} 
We omit the explicit descriptions of the isogenies $\phi,\phi'$, as they are lengthy and will not be used in this paper. 
\bigskip\\
{\bf\underline{Example:}} We will prove the following theorem:
\begin{theorem}\label{exthm}
Suppose that $U=4,V=-4,\Delta=4$, and that $W=\frac{3}{4}-4n$ for some non-negative integer $n$ such that $q:=8n+11,r:=256n^2-2912n-2087$ are both either prime or minus a prime. Then, the Mordell-Weil rank of $J$ and $J'$ over $\Q$ is $2$. 

Let $l$ be a prime such that $2,r$ are inert in $L$, where $L$ is $\Q(\sqrt{l})$. Then, the Mordell-Weil rank of $J$ and $J'$ over $L$ is $4$ if $J$ has a point defined over $L$ such that no multiple of it is defined over $\Q$, and is $2$ otherwise.  
\end{theorem}                
In the theorem, $U,V,W,\Delta$ have been carefully assigned to ensure that $[2]$ is an element of $I$ and is not an element of $I'$, both for $\Q$ and for $L$. The fact that $l$ is congruent to $1$ modulo $4$ ensures that the narrow and wide class numbers of $L$ are the same, which is useful when dealing with the units in $L$. The restrictions on $q,r$ help with computing $I$ and $I'$, although it turns out that we can avoid the restriction that $q$ is inert in $L$.    
\bigskip\\
\pf By substituting directly into the definitions in Theorem~\ref{isogthm}, we see that $C,C'$ are the curves
\begin{eqnarray}  
C\colon Y^2\tabeq(X^2+4X-4)\label{cur1}\\
\tabtimes(4X^4+20X^3-(16n+9)X^2+(32n+14)X-16n-13),\nonumber\\
C'\colon Y^2\tabeq2^{12}\cdot3(X-2)(X-4)\nonumber\\
\tabtimes((80n+58)X^4-(384n+228)X^3+(656n+481)X^2\nonumber\\
\tabminus(480n+438)X+128n+136).\nonumber 
\end{eqnarray}  
By replacing $(x,y)$ by $\lb(\frac{x}{2},\frac{y}{2^8}\rb)$, and altering $\phi,\phi'$ appropriately, we may replace $C'$ by the curve
\begin{eqnarray} 
C'\colon Y^2\tabeq3(X-1)(X-2)\label{cur2}\\
\tabtimes((320n+232)X^4-(768n+456)X^3+(656n+481)X^2\nonumber\\
\tabminus(240n+219)X+32n+34).\nonumber 
\end{eqnarray}   
Adopting notation from the previous section, we define
\begin{eqnarray*}
F_2(X)\tabeq X^2+4X-4,\\
F_4(X)\tabeq4X^4+20X^3-(16n+9)X^2+(32n+14)X-16n-13,\\
F_2'(X)\tabeq(X-1)(X-2),\\
F_4'(X)\tabeq3((320n+232)X^4-(768n+456)X^3+(656n+481)X^2\\
\tabminus(240n+219)X+32n+34).
\end{eqnarray*}
Consider the Mordell-Weil rank over $\Q$. It suffices to show that 
\begin{equation} 
I=\langle[2],[q]\rangle,I'=\langle[r]\rangle;\label{req1} 
\end{equation} 
this implies that the product of the orders of $I$ and $I'$ is $2^3$, which implies that the Mordell-Weil rank is $2$, by Theorem~\ref{prop}. Define $K$ to be $\Q$ (for the purpose of applying the lemmas in Section~\ref{groups}). 

The discriminants of $F_2F_4,F_2'F_4'$ are $2^{17}q^2r^3,2^{19}3^{22}q^3r^2$ respectively, which implies that $J$ and $J'$ have good reduction outside $2,q,r$ (since $J$ and $J'$ are isogenous over $\Q$, their primes of bad reduction are the same). Note that $2,3,q,r$ are mutually distinct, as can be easily shown. By Lemma~\ref{lem1}, we have
\[I\hookrightarrow\langle[-1],[2],[q],[r]\rangle,I'\hookrightarrow\langle[-1],[2],[q],[r]\rangle.\]
To establish the first statement of (\ref{req1}), we need to show that $[2],[q]$ are elements of $I$, and that $[-1],[r],[-r]$ are not. $J$ contains the points $[(2(-1\pm\sqrt2),0)-\infty^+]$, which have the same image under $\phi$, so $[2]$ is an element of $I$. The two points on $J$ which map to the point $[(1,0)-(2,0)]$ under $\phi$ are defined over $\Q(\sqrt{q})$ individually, so $[q]$ is an element of $I$. Lemma~\ref{lem3} implies that $[-1]$ is not an element of $I$, as can easily be seen by inspecting $F_2$ and $F_4$. Lemma~\ref{lem2} with $\mkp=(r)$ implies that $[r],[-r]$ are not elements of $I$; indeed, $F_4$ is not a constant times a square modulo $(r)$, since $r^2$ does not divide the discriminant of $F_4$, which is $2^{12}q^2r$.

To establish the second statement of (\ref{req1}), we need to show that $[r]$ is an element of $I'$, and that $[2],[-2],[2q],[-2q],[q],[-q],[-1]$ are not. The two points on $J'$ which map to the point $[(2(-1+\sqrt2),0)-(2(-1-\sqrt2),0)]$ under $\phi'$ are defined over $\Q(\sqrt{r})$ individually, so $[r]$ is an element of $I'$. Lemma~\ref{lem2} with $\mkp=(2)$ implies that $[2],[-2],[2q],[-2q]$ are not elements of $I'$, as can easily be seen by inspecting $F_2'$ and $F_4'$. Lemma~\ref{lem2} with $\mkp=(q)$ implies that $[q],[-q]$ are not elements of $I'$; indeed, $F_4'$ is not a constant times a square modulo $(q)$, since $q^2$ does not divide the discriminant of $F_4'$, which is $2^53^{10}qr^2$. Finally, the image of $I$ under the natural map $\Q^*/{\Q^*}^2\rightarrow\Q_q^*/{\Q_q^*}^2$ has order $2^2$, since $[2],[q]$ are elements of $I$ and $2,q,2q$ are not squares in $\Q_q$ (by hypothesis, $q$ is congruent to $3$ modulo $8$), so, by Lemma~\ref{lem4} with $\mkp=(q)$, $[D]$ is an element of $I'$ implies that $D$ is a square in $\Q_q$. We deduce that $[-1]$ is not an element of $I'$. 
\bigskip\\
Consider the Mordell-Weil rank over $L$. Let $\mko$ be the ring of integers of $L$, and fix an embedding $\rho\colon L\hookrightarrow\R$. Let $\mu$ be a fundamental unit of $\mko$ such that $\rho(\mu)$ is positive. Define $K$ to be $L$ (for the purpose of applying the lemmas in Section~\ref{groups}).

There are two cases to consider: 
\begin{enumerate}
\item $q$ splits in $L$. 
\item $q$ is inert in $L$.
\end{enumerate}
\underline{Case 1:} Suppose that $q$ splits in $L$. Then, there are prime ideals $\mkq_1,\mkq_2$ of $\mko$ such that $(q)=\mkq_1\mkq_2$. Let $\pi$ be an element of $\mko$ whose norm is a square in $\Q$ multiplied by $q$ and such that $\rho(\pi)$ is positive.    
It suffices to show that
\begin{equation}
\langle[2],[q]\rangle\leq I\leq\langle[2],[q],[\pi]\rangle,I'=\langle[r]\rangle;\label{req2}
\end{equation}
this implies that the product of the orders of $I$ and $I'$ is $2^3$ or $2^4$, which implies that the Mordell-Weil rank is $2$ or $4$, by Theorem~\ref{prop}, and then the statement about the Mordell-Weil rank in Theorem~\ref{exthm} follows.  

By Lemma~\ref{lem1}, we have
\begin{equation}
I\hookrightarrow\langle[-1],[\mu],[2],[q],[\pi],[r]\rangle,I'\hookrightarrow\langle[-1],[\mu],[2],[q],[\pi],[r]\rangle,\label{req3}
\end{equation} 
as we now explain. Indeed, suppose that $[D]$ is an element of $I\cup I'$. Without loss of generality, we may assume that $D$ is an element of $\mko$. Since $2,r$ are inert in $L$, by hypothesis, $(2),(r)$ are prime ideals of $\mko$. So, for some non-negative rational integer $i$, we have
\begin{equation}
(D)=(2)^{s_2}\mkq_1^{s_{q_1}}\mkq_2^{s_{q_2}}(r)^{s_r}\prod_{i}\mkp_i^{t_i},\label{idealfac}
\end{equation} 
where the $\mkp_i$'s are mutually distinct prime ideals of $\mko$ not equal to any of $(2),(q),(r)$, $s_2,s_{q_1},s_{q_2},s_r$ are non-negative rational integers, and the $t_i$'s are positive rational integers. All the $t_i$'s are even; indeed, if some $t_i$ is odd, then $L(\sqrt{D})/L$ is ramified at $\mkp_i$, but $J$ and $J'$ have good reduction at $\mkp_i$, which contradicts Lemma~\ref{lem1}. Equating the norms of both sides of (\ref{idealfac}), we deduce that the norm of $D$ is a square in $\Q$ multiplied by $1,-1,q$ or $-q$, and so that $D$ is a square-free rational integer multiplied by $1,\mu,\pi$ or $\pi\mu$ multiplied by a square in $L$. Without loss of generality, we may assume that $D$ is a square-free rational integer multiplied by $1,\mu,\pi$ or $\pi\mu$. Moreover, without loss of generality, we may assume that $l$ does not divide $D$, since $l$ is a square in $L$. Any rational prime distinct from $l$ is unramified in $L/\Q$, so, if any rational prime distinct from $2,q,r$ divides $D$, then some $t_i$ is odd, which, as we have seen, is impossible. So (\ref{req3}) is established.

To establish the first statement of (\ref{req2}), we need to show that $[2],[q]$ are elements of $I$, and that classes of the form $[(-1)^{s_{-1}}\mu^{s_\mu}\pi^{s_\pi}r^{s_r}]$ ($s_{-1},s_\mu,s_\pi,s_r$ are $0$ or $1$, $s_{-1},s_\mu,s_r$ not all $0$) are not. After an analogous argument to that for the Mordell-Weil rank over $\Q$, taking into account the fact that $\rho(\mu)$ and $\rho(\pi)$ are positive, we are left with needing to show that $[\mu],[\mu\pi]$ are not elements of $I$. Since $C$ is defined over $\Q$, $[\mu]$ is an element of $I$ if and only if its conjugate is. Their product is $[-1]$ (since $\mu$ is a fundamental unit of $\mko$ and $l$ is congruent to $1$ modulo $4$), which is not an element of $I$, so $[\mu]$ is not an element of $I$. Similarly, $[\mu\pi]$ is not.

To establish the second statement of (\ref{req2}), we need to show that $[r]$ is an element of $I'$, and that classes of the form $[(-1)^{s_{-1}}\mu^{s_\mu}2^{s_2}q^{s_q}\pi^{s_\pi}]$ ($s_{-1},s_\mu,s_2,s_q,s_\pi$ are $0$ or $1$, not all $0$) are not. After an analogous argument to that for the Mordell-Weil rank over $\Q$, taking into account the fact that $-1,2,q,2q$ are not squares in $L^c_{\mkq_1}$ (by hypothesis, $q$ is congruent to $3$ modulo $8$), we are left with needing to show that $[\mu],[-\mu]$ are not elements of $I'$. Since $C$ is defined over $\Q$, $[\mu]$ is an element of $I$ if and only if
its conjugate is. Their product is $[-1]$ (since $\mu$ is a fundamental unit of $\mko$ and $l$ is congruent to $1$ modulo $4$), so, if $[\mu]$ is an element of
$I'$, then there is a point on $J'$ defined over $\Q(\sqrt{-1})$ which is sent by $\phi'$ to
a point on $J$ defined over $\Q$, which we established to be impossible when looking at the Mordell-Weil rank over $\Q$. So $[\mu]$ is not an element of $I$. Similarly, $[-\mu]$ is not. 
\bigskip\\
\underline{Case 2:} Suppose that $q$ is inert in $L$. Then, $(q)$ is a prime ideal of $\mko$. 
 
It suffices to show that
\[I=\langle[2],[q]\rangle,\langle[r]\rangle\leq I'\leq\langle[-1],[r]\rangle;\] 
this implies that the product of the orders of $I$ and $I'$ is $2^3$ or $2^4$, which implies that the Mordell-Weil rank is $2$ or $4$, by Theorem~\ref{prop}, and then the statement about the Mordell-Weil $L$-rank follows. The proof of this case is similar (but easier) to the proof of the first case. 
\finpf    
\bigskip\\
For Mordell-Weil ranks over $\Q$, Theorem~\ref{exthm} is easy to apply. For Mordell-Weil ranks over $L$, the major difficulty in applying Theorem~\ref{exthm} is deciding whether or not there is a point on $J$ defined over $L$ but not over $\Q$. Suppose that such a point exists; call it $P$. A multiple of $P$ is defined over $\Q$ if and only if the difference between $P$ and its conjugate is a torsion point, and the group of torsion points on $J$ defined over $L$ can usually be computed easily by applying Proposition~\ref{formal}, so, if no multiple of $P$ is defined over $\Q$, then this can usually be established easily.
\bigskip\\
For a given value of $n$ for which Theorem~\ref{exthm} is applicable for the Mordell-Weil rank over $\Q$, our procedure is to search for rational numbers $m$ such that, defining $x$ to be $\frac{-4m}{2m^2-2m+1}$, $4x^4+20x^3-(16n+9)x^2+(32n+14)x-16n-13$ is minus a square multiplied by a prime number, say $l$; define $L$ to be $\Q(\sqrt{l})$. Then, $x^2+4x-4$ is minus a square, so $J$ has a point defined over $L$ but not over $\Q$, namely $[(x,y)-\infty^+]$, where $y$ is a rational number multiplied by $\sqrt{l}$. Searching for
rational numbers $x$ such that $x^2+4x-4$ is a square does not seem to work very well.       
\bigskip\\
We now discuss how Proposition~\ref{formal} can be usually applied to show that the difference between $[(x,y)-\infty^+]$ and its conjugate, i.e. $[(x,y)-(x,-y)]$, is not a torsion point, which implies that the Mordell-Weil rank of $J$ and $J'$ over $L$ is $4$. More precisely, we will show that the group of torsion points on $J$ defined over $L$ is of order $2$, generated by $[(2(-1+\sqrt2),0)-(2(-1-\sqrt2),0)]$; call this point $P_0$. It suffices to show that the group of torsion points on $J$ defined over $L$ is a $2$-group, since the only point of order $2$ is $P_0$ (as can be easily established), and there is no point of order $4$, since $\lambda'(P)=[r]$ and $r$ is not a
square in $L$. By Proposition~\ref{formal}, this is true if there are two primes $p_1,p_2$ of good reduction of $J$ such that the greatest common divisor of the numbers of points on the reductions of $J$ modulo $p_1,p_2$ is a power of $2$.
 
In the proof of Theorem~\ref{exthm}, we noted that $J$ has good reduction at $3$. It can be easily established that Theorem~\ref{exthm} is applicable for the Mordell-Weil rank over $\Q$ only if $3$ divides $n$, that $3$ splits in $L$ if $3$ divides $n$, and that the number of points on the reduction of $J$ modulo $3$ is $36$ if $3$ divides $n$. So we need to find a prime $p$ of good reduction of $J$, not equal to $3$ and which splits in $L$, such that $3$ does not divide the number of points on the reduction of $J$ modulo $p$.
\bigskip\\
The three smallest values of $n$ for which Theorem~\ref{exthm} is applicable for the Mordell-Weil rank over $\Q$ are $0,6,9$. For each value of $n$, we give $C$, $C'$, $q$, $r$, and several primes $l$ such that the Mordell-Weil rank over $L$ of $J$ and $J'$ is $4$ (more precisely, all the suitable primes of six digits or less, with the reduced numerator and denominator of $m$ lying between $-10$ and $10$). For each prime $l$, we give $m$, a point on $J$ defined over $L$, a prime $p$ chosen as in the previous paragraph, and the number of points on the reduction of $J$ modulo $p$.                    
\begin{center}
${\mathbf{\underline{n=0}}}$ 
\end{center}
\begin{eqnarray*}
C\colon Y^2\tabeq(X^2+4X-4)(4X^4+20X^3-9X^2+14X-13),\\
C'\colon Y^2\tabeq3(X-1)(X-2)(232X^4-456X^3+481X^2-219X+34).
\end{eqnarray*}
\[q=11,r=-2087.\] 
\begin{center}
\begin{tabular}{ccccc} 
$l$&$m$&Point&$p$&$|J(\F_p)|$\\
$13$&$0$&$[(0,2\sqrt{13})-\infty^+]$&$17$&$400$\\
$47269$&$-3$&$\lb[\lb(\frac{12}{25},\frac{238\sqrt{47269}}{15625}\rb)-\infty^+\rb]$&$5$&$62$\\ 
$71341$&$2$&$\lb[\lb(-\frac{8}{5},\frac{14\sqrt{71341}}{125}\rb)-\infty^+\rb]$&$5$&$62$ 
\end{tabular} 
\end{center} 
\begin{center}
${\mathbf{\underline{n=6}}}$
\end{center}
\begin{eqnarray*}
C\colon Y^2\tabeq(X^2+4X-4)(4X^4+20X^3-105X^2+206X-109),\\
C'\colon Y^2\tabeq3(X-1)(X-2)(2152X^4-5064X^3+4417X^2-1659X+226).
\end{eqnarray*}
\[q=59,r=-10343.\]
\begin{center}
\begin{tabular}{ccccc}
$l$&$m$&Point&$p$&$|J(\F_p)|$\\
$658069$&$-\frac{1}{3}$&$\lb[\lb(\frac{12}{17},\frac{14\sqrt{658069}}{4913}\rb)-\infty^+\rb]$&$5$&$62$
\end{tabular}
\end{center}
\begin{center}
${\mathbf{\underline{n=9}}}$
\end{center}
\begin{eqnarray*}
C\colon Y^2\tabeq(X^2+4X-4)(4X^4+20X^3-153X^2+302X-157),\\
C'\colon Y^2\tabeq3(X-1)(X-2)(3112X^4-7368X^3+6385X^2-2379X+322).
\end{eqnarray*}
\[q=83,r=-7559.\]
\begin{center}
\begin{tabular}{ccccc}
$l$&$m$&Point&$p$&$|J(\F_p)|$\\
$157$&$0$&$[(0,2\sqrt{157})-\infty^+]$&$11$&$100$\\ 
$679741$&$2$&$\lb[\lb(-\frac{8}{5},\frac{14\sqrt{679741}}{125}\rb)-\infty^+\rb]$&$5$&$28$  
\end{tabular}
\end{center}

Peter R. Bending,\\  
Institute of Mathematics and Statistics.\\
University of Kent at Canterbury,\\
Canterbury,\\
Kent,\\
England.\\\\
Email: P.R.Bending@ukc.ac.uk 
\end{document}